\documentclass{article}
\usepackage{amssymb}
\usepackage{amsmath}
\usepackage{amsfonts}

\newtheorem{theorem}{Theorem}

\newenvironment{proof}[1][Proof]{\noindent\textbf{#1.} }{\ \rule{0.5em}{0.5em}}
\input{tcilatex}
\begin{document}

\title{\textbf{ON TIMELIKE AND SPACELIKE DEVELOPABLE RULED SURFACES}}
\author{\textbf{Yusuf YAYLI}\thanks{%
Ankara University, Faculty of Science, Department of Mathematics, Ankara,
TURKEY}\textbf{, Semra SARACOGLU}\thanks{%
Siirt University, Faculty of Science and Arts, Department of Mathematics,
Siirt, TURKEY}}
\maketitle

\begin{abstract}
In this study, we have obtained the distribution parameter of a ruled
surface generated by a straight line in Frenet trihedron moving along a
timelike curve and also along another curve with the same parameter. At this
time, the Frenet frames of these timelike curves are not the same. We have
moved the director vector of the first curve along the second curve. It is
shown that the ruled surface is developable if and only if the base timelike
curve is helix. In addition, some smilarities and differences are presented
with theorems and results.

\textbf{AMS Subj. Class.:} 53A04, 53A05, 53B30.

\textbf{Key words:}Timelike curve, distribution parameter, Minkowski space.
\end{abstract}

\section{INTRODUCTION}

In general, it is known that the Lorentz-Minkowski space is the metric space 
$E_{1}^{3}=(%
\mathbb{R}
^{3},\left\langle ,\right\rangle ),$ where the metric $\left\langle
,\right\rangle $ is given by%
\begin{equation*}
\left\langle u,v\right\rangle =u_{1}v_{1}+u_{2}v_{2}-u_{3}v_{3},\text{ \ \ \ 
}u=(u_{1},u_{2},u_{3}),\text{ \ \ \ }v=(v_{1},v_{2},v_{3})
\end{equation*}

The metric $\left\langle ,\right\rangle $ is called the Lorentzian metric
[5]. Also, it can be called $E_{1}^{3}$ as Minkowski space, and $%
\left\langle ,\right\rangle $ as the Minkowski metric.

A surface in 3-dimensional Minkowski space $%
\mathbb{R}
_{1}^{3}=(%
\mathbb{R}
^{3},dx^{2}+dy^{2}-dz^{2})$ is called a timelike surface if the induced
Lorentzian metric on the surface is a non-degenerate. On the other hand, it
can be easily seen that a surface in 3-dimensional Minkowski space $%
\mathbb{R}
_{1}^{3}=(%
\mathbb{R}
^{3},dx^{2}+dy^{2}-dz^{2})$ is called a spacelike surface if the induced
Lorentzian metric on the surface is a positive [1]. A ruled surface is a
surface swept out by a straight line $X$ moving along a curve $\alpha .$ The
various positions of the generating line $X$ are called the rulings of the
surface. Such a surface, thus, has a parametrization in ruled form as
follows,%
\begin{equation}
\Phi (s,v)=\alpha (s)+vX(s).
\end{equation}

We call $\alpha $ to be the base curve, and $X$ to be the director vector.
If the tangent plane is constant along a fixed ruling, then the ruled
surface is called a developable surface. The remaining ruled surfaces are
called skew surfaces[4]. If there exists a common perpendicular to two
preceding rulings in the skew surface, then the foot of the common
perpendicular on the main ruling is called a central point. The locus of the
central points is called the curve of striction.[9]

The ruled surface $M$ is given by the parametrization%
\begin{eqnarray}
\Phi &:&I\times 
\mathbb{R}
\rightarrow 
\mathbb{R}
_{1}^{3} \\
(t,v) &\rightarrow &\Phi (s,v)=\alpha (s)+vX(s)  \notag
\end{eqnarray}%
in $%
\mathbb{R}
_{1}^{3}$ where $\alpha :I\rightarrow 
\mathbb{R}
_{1}^{3}$ is a differentiable timelike curve parametrized by its arc length
in $%
\mathbb{R}
_{1}^{3}$ that is, ($\left\langle \alpha ^{\prime }(s),\alpha ^{\prime
}(s)\right\rangle =-1$) and $X(s)$ is the director vector of the director
curve such that $X$ is orthogonal to the tangent vector field $T$ of the
base curve $\alpha .\left\{ T,N,B\right\} $ is an orthonormal frame field
along $\alpha $ in $%
\mathbb{R}
_{1}^{3}$ where $N$ is the normal vector field of the ruled surface along $%
\alpha $ [10,11].

At this time, with the assistance of $\alpha ,$ we can define curve $\beta .$
Let $\beta (s)$ be a curve and let its parameter be the same as the
parameter of the curve $\alpha (s).$ Nevertheless, we choose the orthonormal
frame field $\left\{ T,N,B\right\} $ along $\beta $ in $IR_{1}^{3}$ where $N$
is the normal vector field of ruled surface along the curve $\beta .$ Thus,%
\begin{equation}
\left\langle T,T\right\rangle =-1,\left\langle N,N\right\rangle
=1,\left\langle X,X\right\rangle =1
\end{equation}%
Smiliarly, as it is said above for the timelike curve $\alpha $, the ruled
surface that produced during the curve $\alpha (s)$is obtained by the
parametrization%
\begin{eqnarray*}
\Phi _{\alpha } &:&I\times 
\mathbb{R}
\rightarrow 
\mathbb{R}
_{1}^{3} \\
(t,v) &\rightarrow &\Phi _{\alpha }(s,v)=\alpha (s)+vX(s)
\end{eqnarray*}%
On the other hand, let $P_{X}$ be distribution parameter of the ruled
surface, then%
\begin{equation}
P_{X}=\frac{\det (\alpha ^{\prime },X,X^{\prime })}{\left\langle X^{\prime
},X^{\prime }\right\rangle }
\end{equation}%
and also, we can get the ruled surface that produced during the curve $\beta
(s)$ with each fixed line $X$ of the moving space $H$ as:%
\begin{equation*}
\Phi _{\beta }(s,v)=\beta (s)+vX(s)
\end{equation*}%
Just after that, the distribution parameter of the ruled surface for the
curve $\beta $ can be given as:

\begin{equation}
\overset{\backsim }{P_{X}}\text{ }=\frac{\det (\beta ^{\prime },X,X^{\prime
})}{\left\langle X^{\prime },X^{\prime }\right\rangle }
\end{equation}

\begin{theorem}
A ruled surface is a developable surface if and only if the distribution
parameter of the ruled surface is zero.[11]
\end{theorem}

\section{BASIC CONCEPTS}

Now we give the basic concepts on differential geometry of cuves and
surfaces in Lorentz-Minkowski space about our study.

\subsection{The Timelike Case [5]}

We suppose that $\alpha $ is a timelike curve. Then $T^{\prime }(s)\neq 0$
is a spacelike vector independent with $T(s).$ We define the curvature of $%
\alpha $ at $s$ as $k_{1}(s)=\left\vert T^{\prime }(s)\right\vert .$The
normal vector $N(s)$ is defined by 
\begin{equation}
N(s)=\frac{T^{\prime }(s)}{\kappa (s)}=\frac{\alpha ^{\prime \prime }(s)}{%
\left\vert \alpha ^{\prime \prime }(s)\right\vert }.
\end{equation}%
Moreover $k_{1}(s)=\left\langle T^{\prime }(s),N(s)\right\rangle .$ We call
the binormal vector $B(s)$ as 
\begin{equation}
B(s)=T(s)\times N(s)
\end{equation}%
The vector $B(s)$ is unitary and spacelike. For each $s,$ $\left\{
T,N,B\right\} $ is an orthonormal base of $E_{1}^{3}$ which is called the
Frenet trihedron of $\alpha .$ We define the torsion of $\alpha .$ We define
the torsion of $\alpha $ at $s$ as%
\begin{equation}
k_{2}(s)=\left\langle N^{\prime }(s),B(s)\right\rangle .
\end{equation}%
By differentiation each one of the vector functions of the Frenet trihedron
and putting in relation with the same Frenet base, we obtain the Frenet
equations, namely,%
\begin{equation}
\left[ 
\begin{array}{c}
T^{\prime } \\ 
N^{\prime } \\ 
B^{\prime }%
\end{array}%
\right] =\left[ 
\begin{array}{ccc}
0 & k_{1} & 0 \\ 
k_{1} & 0 & k_{2} \\ 
0 & -k_{2} & 0%
\end{array}%
\right] \left[ 
\begin{array}{c}
T \\ 
N \\ 
B%
\end{array}%
\right]
\end{equation}

\section{TIMELIKE DEVELOPABLE RULED SURFACES}

Let 
\begin{equation}
\alpha :I\rightarrow 
\mathbb{R}
_{1}^{3}
\end{equation}
be a timelike curve and $\left\{ T,N,B\right\} $ be Frenet vector, where $T,$
$N$ and $B$ are the tangent, principal normal and binormal vectors of the
curve, respectively. $T^{\prime }(s)$ and $B(s)$ vectors are spacelike, and
at the same time $B(s)$ is unitary vector

The two coordinate systems $\left\{ O;T,N,B\right\} $ and $\left\{ O^{\prime
};\overrightarrow{e}_{1},\overrightarrow{e}_{2},\overrightarrow{e}%
_{3}\right\} $ are orthogonal coordinate systems in $%
\mathbb{R}
_{1}^{3}$ which represent the moving space $H$ and the fixed space $%
H^{\prime },$ respectively. Let us express the displacements ($H/H^{\prime }$%
) of $H$ with respect to $H^{\prime }.$ During the one parameter spatial
motion $H/H^{\prime },$ each line $X$ of the moving space $H,$ generates, in
generally, a ruled surface in the fixed space $H^{\prime }.$

On the other hand, let $X$ be a unit vector. Thus%
\begin{equation}
X\in Sp\left\{ T,N,B\right\} \text{ and }X=x_{1}T+x_{2}N+x_{3}B
\end{equation}%
such that 
\begin{equation}
\left\langle X,X\right\rangle =1\text{ }
\end{equation}%
\textbf{\ }We can obtain the distribution parameter of the ruled surface
generated by line $X$ of the moving space $H.$ Different cases can be
investigated as following:

Let the vector $T^{\prime }(s)$ and $B(s)$ be spacelike. Thus,\ from (11) we
have

\begin{equation}
X^{\prime }=x_{1}T^{\prime }+x_{2}N^{\prime }+x_{3}B^{\prime },\text{ \ \ \ }%
-x_{1}^{2}+x_{2}^{2}+x_{3}^{2}=1
\end{equation}%
substituting (10) into (13), then%
\begin{eqnarray}
X^{\prime } &=&x_{1}(k_{1}N)+x_{2}(k_{1}T+k_{2}B)+x_{3}(-k_{2}N) \\
&=&x_{2}k_{1}T+(x_{1}k_{1}-x_{3}k_{2})N+x_{2}k_{2}B  \notag
\end{eqnarray}%
From (4), we obtain%
\begin{equation*}
P_{X}=\frac{\det
(T,x_{1}T+x_{2}N+x_{3}B,x_{2}k_{1}T+(x_{1}k_{1}-x_{3}k_{2})N+x_{2}k_{2}B)}{%
x_{2}^{2}(k_{1}^{2}+k_{2}^{2})+(x_{1}k_{1}-x_{3}k_{2})^{2}}
\end{equation*}%
And then, it can be easily seen that%
\begin{eqnarray}
P_{X} &=&\frac{x_{2}^{2}k_{2}-x_{1}x_{3}k_{1}+x_{3}^{2}k_{2}}{%
x_{2}^{2}(k_{1}^{2}+k_{2}^{2})+(x_{1}k_{1}-x_{3}k_{2})^{2}} \\
&=&\frac{k_{2}(x_{2}^{2}+x_{3}^{2})-x_{1}x_{3}k_{1}}{%
x_{2}^{2}(k_{1}^{2}+k_{2}^{2})+(x_{1}k_{1}+x_{3}k_{2})^{2}}
\end{eqnarray}%
The ruled surface is developable if and only if $P_{X}$ is zero. From (16)%
\begin{equation}
P_{X}=0\text{ if and only if }\frac{k_{1}}{k_{2}}=\frac{x_{2}^{2}+x_{3}^{2}}{%
x_{1}x_{3}}
\end{equation}%
Hence we state the following theorem:

During the one parameter spatial motion $H/H^{\prime }$ the ruled surface in
the fixed space $H^{\prime }$ generated by a line $X$ of the moving space $H$
is developable if and only if $\beta ^{\prime }(s)$ is a helix such that
harmonic curvature $h$ of the base timelike curve $\alpha (s)$ satisfies the
equality 
\begin{equation*}
h=\frac{k_{1}}{k_{2}}=\frac{x_{2}^{2}+x_{3}^{2}}{x_{1}x_{3}}
\end{equation*}%
Different special situations can occur and some results of these special
cases are smiliar in [9,10,11]:

\textbf{1. The Case }$X=T$

In this case $x_{1}=1,$ $x_{2}=x_{3}=0,$ thus from (16)%
\begin{equation*}
P_{X}=0
\end{equation*}%
At this time, we can say that the ruled surface is timelike. According to
all of these the following theorem can be given:

\begin{theorem}
During the one-parameter spatial motion $H/H^{\prime },$the timelike ruled
surface in the fixed space $H^{\prime }$ generated by the tangent line $T$
of the timelike curve $\alpha $ in the moving space $H$ is developable.
\end{theorem}

\textbf{2. The Case }$X=N$

In this case, $x_{2}=1,$ $x_{1}=x_{3}=0,$ thus from (16)%
\begin{equation*}
P_{N}=\dfrac{k_{2}}{k_{1}^{2}+k_{2}^{2}}
\end{equation*}%
If $P_{N}$ is zero then $k_{2}$ is zero. Thus, the timelike curve $\alpha
(s) $ is a planar curve. Hence the following theorem is hold:

\begin{theorem}
During the one-parameter spatial motion $H/H^{\prime },$the timelike ruled
surface in the fixed space $H^{\prime }$ generated by the normal line $N$ of
the timelike curve $\alpha $ in the moving space $H$ is developable.
\end{theorem}

\textbf{3. The Case }$X=B$

In this case, $x_{3}=1,$ $x_{1}=$ $x_{2}=0;$ thus from (16)

\begin{equation*}
P_{N}=\dfrac{k_{2}}{k_{2}^{2}}=\dfrac{1}{k_{2}}
\end{equation*}%
At this time, it can be seen that the ruled surface is timelike.

\textbf{4. The Case }$X$ \textbf{is in the Normal Plane}

In this case, $x_{1}$ is zero. The ruled surface is timelike.and also it is
developable. Now from (16)

\begin{eqnarray*}
P_{X} &=&\dfrac{k_{2}(x_{2}^{2}+x_{3}^{2})}{%
x_{2}^{2}(k_{1}^{2}+k_{2}^{2})+x_{3}^{2}k_{2}^{2}} \\
&=&\dfrac{k_{2}(x_{2}^{2}+x_{3}^{2})}{%
(x_{2}^{2}+x_{3}^{2})k_{2}^{2}+x_{2}^{2}k_{1}^{2}}
\end{eqnarray*}%
If $P_{X}$ is zero hen $k_{2}=0$ $(or$ $x_{2}^{2}+x_{3}^{2}=0).$ Thus, if $%
k_{2}=0,$ then the timelike curve $\alpha (s)$ is a planar curve. Hence the
following theorem can be given as:

\begin{theorem}
During the one-parameter spatial motion $H/H^{\prime },$the timelike ruled
surface in the fixed space $H^{\prime }$ generated by a line $X$ in the
normal plane of $H$ is developable if and only if the timelike curve $\alpha 
$ is a planar curve in normal plane.
\end{theorem}

\textbf{5. The Case }$X$ \textbf{is in the Osculating Plane}

In this case, $x_{3}$ is zero. Thus the ruled surface is timelike and also
it is developable since from (16)%
\begin{eqnarray*}
P_{X} &=&\dfrac{x_{2}^{2}k_{2}}{%
x_{2}^{2}(k_{1}^{2}+k_{2}^{2})+x_{1}^{2}k_{1}^{2}} \\
&=&\dfrac{x_{2}^{2}k_{2}}{(x_{2}^{2}+x_{1}^{2})k_{1}^{2}+x_{2}^{2}k_{2}^{2}}
\end{eqnarray*}%
If $P_{X}=0$ then $X=T.$ This is the \textbf{case.1. }If $k_{2}=0$ then the
timelike curve $\alpha (s)$ is a planar curve. Therefore \textbf{Therorem.3 }%
can be restated.

\textbf{6. The Case }$X$ \textbf{is in the Rectifying Plane}

In this case, $x_{2}$ is zero. From (16), we can give the distribution
parameter of the ruled surface as:

\begin{equation*}
P_{X}=\dfrac{-x_{1}x_{3}k_{1+}x_{3}^{2}k_{2}}{(x_{1}k_{1}-x_{3}k_{2})^{2}}
\end{equation*}%
If $P_{X}=0$ then $\dfrac{k_{1}}{k_{2}}=-\dfrac{x_{3}}{x_{1}}.$ Thus it can
be seen that the timelike curve $\alpha (s)$ is a helix if and only if the
ruled surface is developable such that $\dfrac{k_{1}}{k_{2}}=-\dfrac{x_{3}}{%
x_{1}}.$ At this time, the timelike curve is helix if and only if the base
curve is striction line such that $\dfrac{k_{1}}{k_{2}}=-\dfrac{x_{3}}{x_{1}}%
.$

\subsection{TIMELIKE OR SPACELIKE DEVELOPABLE RULED SURFACES}

As it is said in the first part of the study, with the assistance of $\alpha
,$ another different cuve $\beta $ can be defined with the same parameter of
the timelike curve $\alpha (s),$ such that 
\begin{equation*}
\beta ^{\prime }=\lambda _{1}T+\lambda _{2}N+\lambda _{3}B,
\end{equation*}%
At this time, we can get the ruled surface that produced during the curve $%
\beta (s)$ with each line $X$ of the moving space $H$ as%
\begin{equation*}
\Phi _{\beta }(s,v)=\beta (s)+vX(s)
\end{equation*}%
Smiliarly, the two coordinate systems $\left\{ O;T,N,B\right\} $ and $%
\left\{ O^{\prime };\overrightarrow{e}_{1},\overrightarrow{e}_{2},%
\overrightarrow{e}_{3}\right\} $ are orthogonal coordinate systems in $%
IR_{1}^{3}$ which represent the moving space $H$ and the fixed space $%
H^{\prime },$ respectively. Let us express the displacements ($H/H^{\prime }$%
) of $H$ with respect to $H^{\prime }.$ During the one parameter spatial
motion $H/H^{\prime },$ each line $X$ of the moving space $H,$ generates, in
generally, a ruled surface in the fixed space $H^{\prime }.$ Subsequently,
it should be seen that the curve can be spacelike or timelike but it can not
be null. For this situation, taking $\left\langle \beta ^{\prime
},X\right\rangle =0$ is enough.

Let $X$ be a vector and $X$ is fixed. Thus,

\begin{equation}
X\in Sp\left\{ T,N,B\right\} \text{ and }X=x_{1}T+x_{2}N+x_{3}B
\end{equation}%
such that 
\begin{equation}
\left\langle X,X\right\rangle =\pm 1\text{ }
\end{equation}%
We can obtain the distribution parameter of the ruled surface generated by
line $X$ of the moving space $H.$ As it is investigated in the first kind of
timelike curve $\alpha ,$ let the vector $T^{\prime }(s)$ and $B(s)$ of the
curve $\beta $ be spacelike. Smiliarly from $X=x_{1}T+x_{2}N+x_{3}B,$ we get

\begin{equation}
X^{\prime }=x_{1}T^{\prime }+x_{2}N^{\prime }+x_{3}B^{\prime },\text{ \ \ \ }%
-x_{1}^{2}+x_{2}^{2}+x_{3}^{2}=\pm 1
\end{equation}%
substituting (10) into (20), then

\begin{eqnarray}
X^{\prime } &=&x_{1}T^{\prime }+x_{2}N^{\prime }+x_{3}B^{\prime } \\
&=&x_{1}(k_{1}N)+x_{2}(k_{1}T+k_{2}B)+x_{3}(-k_{2}N)  \notag \\
&=&x_{2}k_{1}T+(x_{1}k_{1}-x_{3}k_{2})N+x_{2}k_{2}B  \notag
\end{eqnarray}%
From (5), we obtain

\begin{eqnarray}
P_{X} &=&\frac{\det (\beta ^{\prime },X,X^{\prime })}{%
x_{2}^{2}k_{1}^{2}+(x_{1}k_{1}+x_{3}k_{2})^{2}+x_{2}^{2}k_{2}^{2}} \\
&=&\frac{\det (\lambda _{1}T+\lambda _{2}N+\lambda
_{3}B,x_{1}T+x_{2}N+x_{3}B,x_{2}k_{1}T+(x_{1}k_{1}-x_{3}k_{2})N+x_{2}k_{2}B)%
}{x_{2}^{2}k_{1}^{2}+(x_{1}k_{1}-x_{3}k_{2})^{2}+x_{2}^{2}k_{2}^{2}} \\
&=&\frac{\lambda _{1}(x_{2}^{2}k_{2}-x_{1}x_{3}k_{1}+x_{3}^{2}k_{2})-\lambda
_{2}(x_{1}x_{2}k_{2}-x_{2}x_{3}k_{1})+\lambda
_{3}(x_{1}^{2}k_{1}-x_{1}x_{3}k_{2}-x_{2}^{2}k_{1})}{%
x_{2}^{2}k_{1}^{2}+(x_{1}k_{1}-x_{3}k_{2})^{2}+x_{2}^{2}k_{2}^{2}}  \notag \\
&=&\frac{\lambda _{1}((x_{2}^{2}+x_{3}^{2})k_{2}-x_{1}x_{3}k_{1})-\lambda
_{2}(x_{1}x_{2}k_{2}-x_{2}x_{3}k_{1})+\lambda
_{3}((x_{1}^{2}-x_{2}^{2})k_{1}-x_{1}x_{3}k_{2})}{%
x_{2}^{2}k_{1}^{2}+(x_{1}k_{1}-x_{3}k_{2})^{2}+x_{2}^{2}k_{2}^{2}}
\end{eqnarray}
\ \ \ 

\subsubsection{\textbf{For }$\protect\lambda _{1}=1,$\textbf{\ }$\protect%
\lambda _{2}=\protect\lambda _{3}=0$ $(\protect\beta ^{\prime }=T)$}

The same special cases and results have occured. We have found the same
results as above investigating the timelike curve $\alpha (s).$

\subsubsection{\textbf{For }$\protect\lambda _{2}=1,$\textbf{\ }$\protect%
\lambda _{1}=\protect\lambda _{3}=0$ $(\protect\beta ^{\prime }=N)$}

\begin{equation}
P_{X}=\frac{-x_{1}x_{2}k_{2}-x_{2}x_{3}k_{1}}{%
x_{2}^{2}k_{1}^{2}+(x_{1}k_{1}-x_{3}k_{2})^{2}+x_{2}^{2}k_{2}^{2}}
\end{equation}%
The ruled surface is developable if and only if $P_{X}$ is zero. From (24)
and (25)

\begin{equation}
P_{X}=0\text{ if and only if }\frac{k_{1}}{k_{2}}=\frac{-x_{1}}{x_{3}}
\end{equation}%
It is the same as the situation \textquotedblleft $T^{\prime }(s)$ is
spacelike and $B(s)$ is timelike\textquotedblright\ in [11]. It can be seen
that the curve is helix if and only if the ruled surface is developable such
that $\dfrac{k_{1}}{k_{2}}=\dfrac{-x_{1}}{x_{3}}.$ And also the curve is
helix if and only if the base curve striction line is constant such that $%
\dfrac{k_{1}}{k_{2}}=\dfrac{-x_{1}}{x_{3}}.$

\textbf{SPECIAL CASES}

\textbf{1. The Case }$X=T$

In this case, $x_{1}=1,x_{2}=x_{3}=0.$Thus $P_{T}=0.$

\begin{theorem}
During the one parameter spatial motion $H/H^{\prime }$ the ruled surface in
the fixed space $H^{\prime }$ generated by the tangent line $T$ of the curve 
$\beta (s)$ in the moving space $H$ is developable.
\end{theorem}

\textbf{2. The Case }$X=N$

In this case, $x_{2}=1,x_{1}=x_{3}=0.$Thus $P_{N}=0.$

\begin{theorem}
During the one parameter spatial motion $H/H^{\prime }$ the ruled surface in
the fixed space $H^{\prime }$ generated by the normal line $N$ of the curve $%
\beta (s)$ in the moving space $H$ is developable.
\end{theorem}

\textbf{3. The Case }$X=B$

In this case, $x_{3}=1,x_{1}=x_{2}=0.$Thus $P_{B}=0.$

\begin{theorem}
During the one parameter spatial motion $H/H^{\prime }$ the ruled surface in
the fixed space $H^{\prime }$ generated by the binormal line $B$ of the
curve $\beta (s)$ in the moving space $H$ is developable.
\end{theorem}

\textbf{4. The Case }$X$ \textbf{is in the Normal Plane}

In this case, $x_{1}=0.$The ruled surface is developable, from (25),%
\begin{equation*}
P_{X}=\frac{-x_{2}x_{3}k_{1}}{%
x_{2}^{2}(k_{1}^{2}+k_{2}^{2})+x_{3}^{2}k_{2}^{2}}
\end{equation*}%
If $P_{X}=0$ then $k_{1}=0.$ So, the base curve $\beta (s)$ is a line. The
ruled surface is developable if and only if $x_{2}=0$ or $x_{3}=0.$

\textbf{5.The Case }$X$ \textbf{is in the Osculating Plane}

In this case, \textbf{\ }$x_{3}=0.$The ruled surface is developable, since
from (25),%
\begin{equation*}
P_{X}=\frac{-x_{1}x_{2}k_{2}}{%
x_{2}^{2}(k_{1}^{2}+k_{2}^{2})+x_{1}^{2}k_{1}^{2}}
\end{equation*}%
If $P_{X}$ is zero then $k_{2}=0.$ Thus if $k_{2}=0$ then $\beta (s)$ is a
planar curve. Hence following theorem in the first case can be occured:

\begin{theorem}
During the one parameter spatial motion he ruled surface in the fixed space $%
H^{\prime }$ generated by a line $X$ in the normal plane of $H$ is
developable if and only if the curve $\alpha (s)$ is a planar curve in
osculating plane.
\end{theorem}

\textbf{6.The Case }$X$ \textbf{is in the Rectifying Plane}

In this case, $x_{2}=0,$ from (25) the distribution parameter of the ruled
surface is: 
\begin{equation*}
P_{X}=0
\end{equation*}%
If $P_{X}$ is zero then $X=T.$ This is \textbf{the case 1.}

\subsubsection{\textbf{For }$\protect\lambda _{3}=1,$\textbf{\ }$\protect%
\lambda _{1}=\protect\lambda _{2}=0$ $(\protect\beta ^{\prime }=B)$}

\begin{equation}
P_{X}=\frac{k_{1}(x_{1}^{2}-x_{2}^{2})-x_{1}x_{3}k_{2}}{%
x_{2}^{2}(k_{1}^{2}+k_{2}^{2})+(x_{1}k_{1}-x_{3}k_{2})^{2}}
\end{equation}%
The ruled surface is developable if and only if $P_{X}$ is zero. As a result
we can easily say that

\textbf{Result.1} $P_{X}=0$ if and only if $%
k_{1}(x_{1}^{2}-x_{2}^{2})=x_{1}x_{3}k_{2}$

\textbf{Result.2 }$P_{X}=0$ if and only if $\dfrac{k_{1}}{k_{2}}=\dfrac{%
x_{1}x_{3}}{x_{1}^{2}-x_{2}^{2}}$ and also it can be seen that the timelike
curve $\alpha $ is helix.

\textbf{SPECIAL CASES}

\textbf{1. The Case }$X=T,$

In this case, $x_{1}=1,x_{2}=x_{3}=0.$Thus, from (27),%
\begin{equation*}
P_{T}=\dfrac{k_{1}}{x_{1}^{2}k_{1}^{2}}=.\dfrac{1}{k_{1}}
\end{equation*}%
It can be seen that the ruled surface is not developable.

\textbf{2. The Case }$X=N,$

In this case, $x_{2}=1,x_{1}=x_{3}=0.$So, from (27),

\begin{equation*}
P_{N}=\dfrac{-k_{1}}{k_{1}^{2}+k_{2}^{2}}
\end{equation*}%
If $P_{N}$ is zero then $k_{1}$ is zero.\textbf{\ }In this case, that can
not happen. For this reason, the ruled surface is not developable. And also,
it can be seen that the timelike curve $\alpha $ is Mannheim curve if and
only if $P_{N}$ is constant.

\textbf{3. The Case }$X=B,$

In this case, $x_{3}=1,x_{1}=x_{2}=0.$Thus, from (27),%
\begin{equation*}
P_{B}=0
\end{equation*}

\begin{theorem}
During the one parameter spatial motion $H/H^{\prime }$ the ruled surface in
the fixed space $H^{\prime }$ generated by the binormal line $B$ of the
curve $\beta (s)$ in the moving space $H$ is developable.
\end{theorem}

\textbf{4.The Case }$X$ \textbf{is in the Normal Plane}

In this case, \textbf{\ }$x_{1}=0,$ so%
\begin{equation*}
P_{X}=\dfrac{-k_{1}x_{2}^{2}}{%
x_{2}^{2}(k_{1}^{2}+k_{2}^{2})+x_{3}^{2}k_{2}^{2}}
\end{equation*}%
If $P_{X}=0,$ then $-k_{1}x_{2}^{2}=0.$From here, $x_{2}$ is zero. Hence,
this case gives us \textbf{The Case }$X=B.$

\textbf{5.The Case }$X$ \textbf{is in the Osculating Plane}

In this case, \textbf{\ }$x_{3}=0.$ Then 
\begin{equation*}
P_{X}=\dfrac{k_{1}(x_{1}^{2}-x_{2}^{2})}{%
x_{2}^{2}(k_{1}^{2}+k_{2}^{2})+x_{1}^{2}k_{1}^{2}}
\end{equation*}%
$P_{X}=0$ if and only if $x_{1}=\pm x_{2}.$ From here, it can be easily seen
that if $x_{1}=\left\vert x_{2}\right\vert ,$ then the ruled surface is
developable.

\textbf{6.The Case }$X$ \textbf{is in the Rectifying Plane}

In this case,\textbf{\ }$x_{2}=0.$%
\begin{eqnarray*}
P_{X} &=&\dfrac{k_{1}x_{1}^{2}-x_{1}x_{3}k_{2}}{(x_{1}k_{1}-x_{3}k_{2})^{2}}
\\
&=&\frac{x_{1}(x_{1}k_{1}-x_{3}k_{2})}{(x_{1}k_{1}-x_{3}k_{2})^{2}} \\
&=&\frac{x_{1}}{(x_{1}k_{1}-x_{3}k_{2})}
\end{eqnarray*}%
If $P_{X}=0$ if and only if $x_{1}=0.$ Then $X=x_{3}B$ and $x_{3}^{2}=1.$ At
this time, the ruled surface is developable.

\section{THE CURVE OF THE STRICTION OF A RULED SURFACE}

In this section, as it is said in the first part of the study, the various
positions of the generating line $X$ are called the rulings of the surface.
Such a surface, thus, has a parametrization in ruled form as follows,%
\begin{equation}
\Phi _{\alpha }(s,v)=\alpha (s)+vX(s).
\end{equation}%
the curve of the striction of a ruled surface can be given. This can be
investigated for the base timelike curve $\alpha $ and also for another
different curve $\beta $. If there exists a common perpendicular to two
preceding rulings in the ruled surface, then the foot of the common
perpendicular on the main ruling is called a central point. The locus of the
central points is called the curve of striction.[9]

The ruled surface is obtained by the parametrization%
\begin{eqnarray}
\Phi _{\alpha } &:&I\times 
\mathbb{R}
\rightarrow 
\mathbb{R}
_{1}^{3} \\
(t,v) &\rightarrow &\Phi _{\alpha }(s,v)=\alpha (s)+vX(s)  \notag
\end{eqnarray}%
On the other hand, we can get the ruled surface that produced during the
curve $\beta (s)$ with each line $X$ of the moving space $H$ as:%
\begin{equation}
\Phi _{\beta }(s,v)=\beta (s)+vX(s)
\end{equation}%
Thus, the curve of the striction of a ruled surface can be given by:%
\begin{equation}
\widetilde{\beta }(s)=\beta (s)-\frac{\left\langle \dfrac{d\beta }{ds}%
,X^{\prime }\right\rangle }{\left\langle X^{\prime },X^{\prime
}\right\rangle }.X
\end{equation}%
In this case, some important situations occur, such as it can be easily seen
that:%
\begin{equation}
\left\langle \dfrac{d\beta }{ds},X^{\prime }\right\rangle =0
\end{equation}%
That is,%
\begin{equation}
\left\langle \lambda _{1}T+\lambda _{2}N+\lambda
_{3}B,x_{2}k_{1}T+(x_{1}k_{1}-x_{3}k_{2})N+x_{2}k_{2}B\right\rangle =0
\end{equation}%
and then,%
\begin{equation}
-\lambda _{1}x_{2}k_{1}+\lambda _{2}(x_{1}k_{1}-x_{3}k_{2})+\lambda
_{3}x_{2}k_{2}=0
\end{equation}%
Accordingly, the following theorem is obtained:

\begin{theorem}
$\widetilde{\beta }=\beta $ if and only if $\left\langle \dfrac{d\beta }{ds}%
,X^{\prime }\right\rangle =0$.
\end{theorem}

\textbf{Result.1 }If the director vector $X$ and the curve $\beta ^{\prime }$
are in the rectifying plane of the timelike curve $\alpha ,$ then $%
\widetilde{\beta }=\beta .$

\begin{proof}
If $X$ is in the rectifying plane, then $x_{2}=0.$ In addition to this, if $%
\beta ^{\prime }$ is in the rectifying plane, then $\lambda _{2}=0.$
\end{proof}

\textbf{Result.2 }If the director vector $X$ \ is in rectifying plane and
the timelike curve of $\alpha $ is a helix, such that $\dfrac{k_{1}}{k_{2}}=%
\dfrac{x_{3}}{x_{1}},$ then the ruled surface is cylinder. Because,%
\begin{equation}
X^{\prime }=(x_{1}k_{1}-x_{3}k_{2})N
\end{equation}%
and%
\begin{equation}
X=x_{1}T+x_{3}B.
\end{equation}%
Thus,%
\begin{equation}
X^{\prime }=0
\end{equation}%
then the director vector $X$ is constant.

\textbf{Result.3.} On the other hand, if the director vector $X$ is in
rectifying plane and the timelike curve $\alpha $ is a helix such that the
curvature $\dfrac{k_{1}}{k_{2}}=\dfrac{x_{3}}{x_{1}},$ then the curve $\beta
(s)$ is the striction curve of the ruled surface.

\section{CONCLUSIONS}

In this paper, the distribution parameter of a ruled surface has been
investigated by taking a timelike curve and also another curve with the same
parameter, with their orthonormal frame fields. By moving the director
vector of the first curve along the second curve, some similarities and
differences are found. Depending on this, important results and theorems are
presented about the cases of the base timelike curve. It is seen that,
according to the timelike curve, the values of distribution parameter of
ruled surfaces are changed.

\textbf{References}

\begin{itemize}
\item[{\textbf{[1]}}] J.K.Beem and P.E. Ehrlich, Global Lorentzian Geometry,
Marcel Dekker Inc. New York, 1981.

\item[{\textbf{[2]}}] Hac\i saliho\u{g}lu, H. H., Differential Geometri,
Ankara \"{U}niversitesi, Fen Fak\"{u}ltesi, 1993.

\item[{\textbf{[3]}}] Hac\i saliho\u{g}lu, H. H., Differential Geometri-II, A 
\"{U}., Fen Fak., 1994.

\item[{\textbf{[4]}}] Hac\i saliho\u{g}lu, H. H., Turgut, A., On the
Distribution Parameter of Spacelike Ruled Surfaces in the Minkowski 3-space,
Far East J.Math Sci 5, 321-328, 1997.

\item[{\textbf{[5]}}] Lopez, R., Differential Geometry of Curves and Surfaces
in Lorentz-Minkowski Space,Preprint: 2008, arXiv:0810.3351v1[math.DG.]

\item[{\textbf{[6]}}] \"{O}zy\i lmaz, E., Yayl\i , Y., On the Closed Motions
and Closed Spacelike Ruled Surfaces, Commun. Fac. Sci. Univ. Ank. Ser. A1
Math. Stat. 49(2000), no. 1-2, 49-58.

\item[{\textbf{[7]}}] T. Ikawa, On curves and submanifolds in an
indefinnite-Riemannian Manifold, Tsukaba j.math 9, 353-371, 1985.

\item[{[8]}] Yayl\i , Y., Saracoglu, S., On Developable Ruled Surfaces,
Algebras, Groups and Geometries (accepted).

\item[{[9]}] Yayl\i , Y., On The Motion of the Frenet Vectors and Spacelike
Ruled Surfaces in the Minkowski 3-Space, Mathematical\&Computational
Applications, Vol.5(1), pp. 49-55, 2000.

\item[{[10]}] Yayl\i , Y., On The Motion of the Frenet Vectors and Timelike
Ruled Surfaces in the Minkowski 3-Space, Dumlup\i nar \"{U}niversitesi Fen
Bilimleri Dergisi, ISSN:1302-3055, Say\i -1, 1999.

\item[{[11]}] Yayl\i , Y., Saracoglu, S., On Developable Ruled Surfaces in
Minkowski Spaces, Advances in Applied Clifford Algebras, DOI
10.1007/s00006-011-0305-5, 2011
\end{itemize}

\end{document}